\documentstyle{amsppt}
\magnification\magstep1
\NoRunningHeads
\TagsOnRight
\vsize 23 true cm
\hsize 15 true cm
\hcorrection{0.6 true cm}

\centerline{\bf ON THE CONJUGACY SEPARABILITY}
\centerline{\bf IN THE CLASS OF FINITE $P$-GROUPS}
\centerline{\bf OF FINITELY GENERATED NILPOTENT GROUPS}
\bigskip

\centerline{\bf E.~A.~Ivanova}
\medskip
\centerline
{\sl Department of Mathematics, Ivanovo State University}
\centerline
{\sl Ermaka St. 39, Ivanovo, 153025, Russia}
\bigskip

\hangindent=28pt
\hangafter=-6
{\it It is proved that for any prime $p$ a finitely generated nilpotent group is conjugacy separable
in the class of finite $p$-groups if and only if the torsion subgroup
of it is a finite $p$-group and the quotient group by the torsion
subgroup is abelian.}
\bigskip

1. Let $\Cal K$ be a class of groups. A group $G$ is called
residual $\Cal K$ (or $\Cal K$-residual) if for each non-unit
element
$a\in G$ there is a homomorphism $\varphi$ of $G$ onto some group $X$ from the
class $\Cal K$ (or $\Cal K$-group) such that  $a\varphi$ is not
the identity.
The group $G$ is called
conjugacy separable in the class $\Cal K$ (or conjugacy $\Cal K$-separable),
if whenever $a$ and $b$ are not conjugate in $G$,
there is a homomorphism $\varphi$ of $G$ onto $\Cal K$-group $X$ such
that $a\varphi$ and $b\varphi$ are not conjugate in $X$.

It is easy to see that a conjugacy $\Cal K$-separable group is also
$\Cal K$-residual. Since in general the inverse statement
is not true, it is interesting to find such classes of
groups for which the property to be $\Cal K$-residual implies conjugacy $\Cal
K$-separability.  The most investigated (and chronologically
first) is the case when class $\Cal K$ is the class $\Cal
F$ of all finite groups; in this case one studies finite residuality and
conjugacy separability respectively. K.~Gruenberg [4] showed that
finitely generated
nilpotent groups are residually finite, and then N.~Blackburn [5] proved that
such groups are conjugacy separable. On the other hand the famous theorem of P.~Hall
states that every finitely generated metabelian group is residually finite,
but there exists  the example of finitely generated metabelian group that is
not conjugacy separable constructed by M.~I.~Kargapolov and E.~I.~Timoshenko [1].

When a class of groups is proved to be $\Cal F$-residual or conjugacy $\Cal
F$-separable, the question arises if these groups are
$\Cal K$-residual or conjugacy $\Cal K$-separable for some subclass $\Cal K$
of the class $\Cal F$. From this point of view the class $\Cal F_p$ of all
finite $p$-groups is frequently considered. For example, A.~I.~Mal'cev [2] showed that
free groups are residually finite. K.~Gruenberg [4] proved that for any
prime $p$ every finitely generated torsion-free nilpotent group is $\Cal
F_p$-residual. From this assertion and from the theorem of W.~Magnus about
$\Cal N$-residuality of free groups (where $\Cal N$ is the class of all
finitely generated torsion-free nilpotent groups) follows that
every free group is $\Cal F_p$-residual for all primes $p$. Since [3] these
groups are conjugacy $\Cal F_p$-separable .

As we have noted above any finitely generated torsion-free nilpotent group is $\Cal
F_p$-residual for every prime $p$. Nevertheless, here we'll show that for every prime $p$
a finitely generated torsion-free nilpotent group is conjugacy
$\Cal F_p$-separable if and only of it is an abelian group. More exactly, we
establish the following statement:

\proclaim{\indent Theorem} Suppose $p$ is a prime.
A finitely generated nilpotent group $G$ is conjugacy $\Cal F_p$-separable
if and only if its torsion subgroup $\tau (G)$ is a $p$-group and the
quotient group $G/\tau (G)$ is abelian.
\endproclaim

\bigskip

2. We'll begin the proof of this theorem from one general remark.

A subgroup $H$ of a group $G$ is called
$\Cal K$-separable if for every $a\in G$, $a\notin H$,
there is a homomorphism $\varphi$ of $G$ onto $\Cal K$-group such
that $a\varphi\notin H\varphi$. It is well known that if
a class $\Cal K$ is closed under homomorphic images, then for every
normal subgroup $N$ of the group $G$ the quotient group $G/N$ is $\Cal
K$-residual if and only if the subgroup $N$ is $\Cal K$-separable.
Consequently, if a class $\Cal K$ is closed under homomorphic images,
subgroups and direct products of finite number of factors, then
a quotient group of the $\Cal K$-residual group by any finite normal subgroup
is $\Cal K$-residual.

To receive the similar result
for the property to be conjugacy $\Cal K$-separable we'll introduce the
following notion.

A subset $M$ of a group $G$ is called conjugacy $\Cal K$-separable in the
group $G$ if
whenever $a\in G$ is not conjugate to any element from $M$, there
is a homomorphism of $G$ onto $\Cal K$-group $X$ such that
$a\varphi$ is not conjugate to any element from $M\varphi$.
It is obvious that subset $M$ is conjugacy $\Cal K$-separable in $G$ if and
only if for every element $a\in G$ which is not conjugate to any element from $M$
there exists a normal subgroup $H$ in $G$ that quotient group
$G/H$ is $\Cal K$-group and $a$ is not conjugate to any element from $MH$.
Evidently also that a group $G$ is conjugacy $\Cal K$-separable iff every
one-element subset of $G$ is conjugacy $\Cal K$-separable in $G$.

\proclaim{\indent Proposition 1} Suppose $\Cal K$ is a class of groups which
is closed under homomorphic images. For every group $G$ and normal subgroup $N$
of $G$ the quotient group $G/N$ is conjugacy $\Cal K$-separable if and only
if every coset of $G$ modulo $N$ is conjugacy $\Cal K$-separable in $G$.
\endproclaim

\demo{\indent Proof} Firstly we'll show  that if some
coset of $G$ modulo $N$ is not conjugacy $\Cal K$-separable in $G$,
then the quotient group $G/N$ is not conjugacy $\Cal K$-separable.

If $a\in G$ is not conjugate to any element from the
coset $bN$, then $aN$ and $bN$  are not conjugate in the
quotient group $G/N$. On the other hand
suppose that for any homomorphism from $G$ onto some $\Cal K$-group the image of $a$
is conjugate to image of some element from a coset $bN$. Let now
$\varphi$ be a homomorphism from $G/N$ onto some
$\Cal K$-group $X$. Then the composition of the natural homomorphism
$\varepsilon: G\to G/N$ and $\varphi$ maps $G$ onto
$X$. Consequently, for some elements $g\in G$ and $c\in N$ we have
$$
(g(\varepsilon\varphi))^{-1}a(\varepsilon\varphi)g(\varepsilon\varphi)=
(bc)(\varepsilon\varphi),
$$
i.~e.
$$
((gN)\varphi)^{-1}(aN)\varphi(gN)\varphi=(bN)\varphi.
$$
Hence, the images of $aN$ and $bN$ under any homomorphism from $G/N$ onto
$\Cal K$-group are conjugate.

Conversely, suppose any coset of $G$ modulo $N$
is conjugacy $\Cal K$-separable in $G$. If elements $aN$ and $bN$
of the quotient group $G/N$ are not conjugate in it, then $a$ is not conjugate
to any element from the coset $bN$ in $G$. Consequently,
$a$ is not conjugate in $G$ to any element from the coset $bNH$ for some normal
subgroup $H$ of the group $G$ such that the quotient group
$G/H$ is $\Cal K$-group. It means that the elements $aNH$ and
$bNH$ of the homomorphic image $G/NH$ of $G/N$ are not conjugate.
Since the class $\Cal K$ is closed under homomorphic images, the quotient
group $G/NH\simeq (G/H)/(NH/H)$ is $\Cal K$-group and the proposition 1 is
proved.
\enddemo

It follows immediately  from the proposition 1 that

\proclaim{\indent Proposition 2} Suppose $\Cal K$ is a class of groups which
is closed under homomorphic images, subgroups and direct products of the finite
number of factors. If a group $G$ is conjugacy $\Cal K$-separable, then a
quotient group $G/N$ of $G$ by any finite normal subgroup $N$ of $G$ is
conjugacy $\Cal K$-separable.
\endproclaim

Actually, it follows from the conditions of the proposition 2 that in any group
$G$ the family of all normal subgroups of $G$ such that the quotient groups by them
are $\Cal K$-groups is closed under finite intersections. Hence, if group $G$
is conjugacy $\Cal K$-separable, then every finite subset of $G$ is
conjugacy $\Cal K$-separable in $G$. So every coset of  $G$ modulo
normal finite subgroup $N$ of $G$ is conjugacy $\Cal K$-separable in $G$.
\smallskip

Now, let $p$ be a prime and $G$ be a finitely generated nilpotent group which
is conjugacy $\Cal F_p$-separable. Since in this case $G$ is $\Cal
F_p$-residual, its torsion subgroup $\tau (G)$ must be a $p$-group.
The quotient group $G/\tau (G)$ is torsion-free and by the
proposition 2 must be conjugacy $\Cal F_p$-separable. So to show that
necessary condition of the theorem is true it is enough to prove

\proclaim{\indent Proposition 3} Suppose $G$ is a finitely generated torsion-free
nilpotent group. If for some prime $p$ the group $G$ is conjugacy $\Cal
F_p$-separable, then $G$ is an abelian group.
\endproclaim

\demo{\indent Proof} Let
$1=Z_0\leqslant Z_1\leqslant \cdots \leqslant Z_r=G$
be the upper central series of $G$. If, on the contrary, $G$ is not abelian,
then $r\geqslant 2$. Let $a$ be an element from $Z_2$ but not from $Z_1$.
Since for any element $g\in G$ the commutator $[a, g]$ is in $Z_1$,
we have $g^{-1}ag=az$ for some $z\in Z_1$. Since the element $a$ is not in
the centre $Z_1$ of the group $G$, there is an element $b\in G$ such
that $b^{-1}ab=ac$ for some non-unit element $c\in Z_1$.

Suppose $q$ is a prime which is not equal to $p$.
Since $Z_1$ is a free abelian group and $c\in Z_1$ is not equal to 1,
for some integer $n\geqslant 1$ the equation $x^{q^n}=c$ has no solutions
in $Z_1$. We assert that the elements $a^{q^n}$ and $a^{q^n}c$ are not conjugate
in $G$.

Indeed, suppose that for some $g\in G$ the equality $g^{-1}a^{q^n}g=a^{q^n}c$
is valid. As $g^{-1}ag=az$ for some $z\in Z_1$, we have
$$
a^{q^n}c=(g^{-1}ag)^{q^n}=a^{q^n}z^{q^n}.
$$
Thus, $c=z^{q^n}$, but this is impossible.

On the other hand let's show that the images of the elements $a^{q^n}$ and
$a^{q^n}c$ under every homomorphism from $G$ onto finite $p$-group are
conjugate.

Suppose $N$ is a normal subgroup of $G$ such that $c^{p^m}\equiv 1\pmod N$
for some integer $m\geqslant 0$. Since numbers $q^n$ and $p^m$ are coprime,
there exists integer $k$ such that $q^nk\equiv 1 \pmod {p^m}$.
As $c$ is a central element, from the equality
$b^{-1}ab=ac$ follows that $b^{-k}ab^k=ac^k$, and therefore
$$
a^{q^n}c\equiv a^{q^n}c^{q^nk}=(ac^k)^{q^n}=b^{-k}a^{q^n}b^k \pmod N.
$$

Thus, the group $G$ is not conjugacy $\Cal F_p$-separable and the proposition
3 is proved.
\enddemo

Note that if torsion subgroup $\tau (G)$ of a finitely generated
nilpotent group $G$ is a finite $p$-group, then from [4] the group $G$ is
$\Cal F_p$-residual. So the inverse statement of the theorem is contained
in the following more general result:

\proclaim{\indent Proposition 4} Suppose a group $G$ has a finite normal
subgroup $F$ such
that the quotient group $G/F$ is a finitely generated abelian group. If
for some prime $p$ the group $G$ is $\Cal F_p$-residual, then $G$ is
conjugacy $\Cal F_p$-separable.
\endproclaim

\demo{\indent Proof}
Let $a$ and $b$ be non-conjugate elements of $G$. We'll show that there exists
a homomorphism $\varphi$ of $G$ onto some finite $p$-group such that the images
$a\varphi$ and $b\varphi$ are not conjugate in $G\varphi$.

By the remark above the quotient group $G/F$ is $\Cal F_p$-residual.
Consequently, if $aF\neq bF$, then the existence of required homomorphism
is evident.

Suppose then that $aF=bF$, i.~e. $b=af$ for some $f\in F$.
Since $G$ is $\Cal F_p$-residual,
there is a normal subgroup $N$ of $G$ that $G/N$ is finite $p$-group and
$N\cap F=1$. We claim that the elements $aN$ and $bN$ of the quotient group
$G/N$ are not conjugate in it.

Indeed, suppose on the contrary that for some element $g\in G$ the equality
$g^{-1}ag=b$ modulo $N$ is valid. As the quotient group
$G/F$ is abelian, we have $g^{-1}ag=ax$ for some $x\in F$.
Then $x=f$ modulo $N$, and as $N\cap F=1$, we get
$x=f$. But it means that the elements $a$ and $b$ are conjugate
in $G$.

The proposition 4 is proved and so the theorem is.
\enddemo

\bigskip
\centerline{\bf References}
\medskip

\noindent\item{1.}
{\it Kargapolov M.~I., Timoshenko E.~I.} On the question of the conjugacy
separability of metabelian groups, IV-th All-Union symposium on group theory
(Novosibirsk, 5 - 9 February, 1973), Lecture notes, Novosibirsk, 1973,  P.~86 -- 88. (Russian)

\noindent\item{2.}
{\it Mal'cev A.~I.} About isomorphic presentation of infinite groups by
matrices, Math. sbornik, V. 8, 1940, P.~405-422. (Russian)

\noindent\item{3.}
{\it Remeslennikov V.~N.} Conjugacy separability of groups, Siberian Math. J.,
V. 12, 1971, P.~1085-1099. (Russian)

\noindent\item{4.}
{\it Gruenberg K.~W.} Residual properties of infinite soluble groups,
Proc. London Math. Soc., V. 7, 1957, P.~29--62.

\noindent\item{5.}
{\it Blackburn N.} Conjugacy in nilpotent groups, Proc. Amer. Math. Soc.,
V. 16, 1965, P.~143--148.
\end